\documentclass[a4paper,enabledeprecatedfontcommands]{scrartcl}
\usepackage{colortbl,amsmath,amsthm,amsfonts,graphicx}

\newtheorem{thm}{Theorem}[section]

\theoremstyle{definition}

\theoremstyle{remark}

\numberwithin{equation}{section}

\begin{document}

\newcommand\Emag{E^{\mathrm{mag}}}
\newcommand\Evoll{E}

\title {Homogenization in magnetic-shape-memory polymer composites}

\author{Sergio Conti\thanks{University of Bonn} \and Martin Lenz\footnotemark[1] \and Matth\"aus Pawelczyk\thanks{TU Dresden} \and Martin Rumpf\footnotemark[1]}

\maketitle

\begin{abstract}  { Magnetic-shape-memory materials (e.g. specific NiMnGa alloys)
react with a large change of shape to the presence of an external magnetic field.
 As an alternative for the difficult to manifacture single crystal of these alloys}
we study composite materials in which small magnetic-shape-memory particles are 
embedded in a polymer matrix. 
The macroscopic properties of the composite depend
strongly on the geometry of the microstructure and on the characteristics of the particles
and the polymer.

We present a variational model based on micromagnetism and elasticity, and derive via homogenization
an effective macroscopic model under the assumption that the microstructure is periodic. We then study numerically the
resulting cell problem, and discuss the effect of the microstructure on the macroscopic material behavior. Our results may be used to 
optimize the shape of the particles and the microstructure.
\end{abstract}

\section{Introduction}
\label{secintro}

The peculiar mechanical behavior of magnetic shape memory (MSM) alloys renders them interesting for many applications. At the same time, they give rise to a new class of models, whose mathematical study is challenging. MSM alloys  are multiferroic materials, which combine a magnetic phase transition with a shape-memory one. The deformation driven by a magnetic field deformation is due to a coupling between the two order parameters, which arises from the difference in magnetic anisotropy between the elastic variants. A typical example is a NiMnGa alloy \cite{Ullakko1996,Tickle99,Murray2000,Sozinov2002,heczko2009magnetic}.

Shape-memory metals are characterized by a  phase transformation from a high-tem\-perature,
high-symmetry phase, called austenite, to a low-tempera\-ture, low-symmetry phase, called martensite. 
Since the point group of the martensite lattice is smaller than in the austenite phase, there are a number of different variants of the martensitic phase, which are obtained from the austenite via different but symmetry-related eigenstrains. 
Whereas in ceramic materials elastic phase transformations typically generate spontaneous strains of less than one percent, in alloys the deformations can be significant. For example in NiMnGa single crystals the deformation between two martensitic variants is about $11\%$.
The practical exploitation of the shape-memory effect for actuation has been hindered by the fact that the transition between austenite and martensite can often be induced only by slow heating and cooling processes. 
 {Furthermore,  it is difficult to select one variant over the other. This leads in typical situations to low-temperature states characterized by a fine mixture of the different variants, in which
martensite appears without significant macroscopic deformation.}

The coupling of the structural phase transformation to a magnetic field renders the materials much easier to control. 
The difference in anisotropy between the different variants of the martensitic phase leads to a
coupling between the magnetization and the eigenstrain and provides a simple mechanism for selecting one of the variants over the others.
At least in clean single crystals, large strains up to $11\%$ have been induced by the application 
of external magnetic fields  \cite{Ullakko1996,Tickle99,Murray2000,Sozinov2002}.

It is however practically very difficult to produce single crystals of good quality. 
Additionally, single crystals tend to be brittle. At the same time, in polycrystals the deformation is blocked by incompatibilities at the grain boundaries. 
It has therefore been proposed 
to use small single-crystal MSM particles embedded in a softer matrix \cite{Feuchtwanger2003,Hosoda2004,Feuchtwanger2005,Scheerbaum2007b,TianChen2009}, 
although magnetically induced strains are typically {smaller} than in single crystals \cite{KauffmannScheerbaumLiuetal2012,LiuScheerbaumKauffmannetal2012,TiChTo14}.
The composite geometry opens the way for optimization of the particle shape and in general of the microstructure.

In this paper  we review recent mathematical progress and extend the analysis of MSM-polymer composites.
A model based on micromagnetism and elasticity was developed in \cite{ContiLenzRumpf2007}, and will be presented in Section \ref{secmodel} below. 
Under the assumption of periodicity of the microstructure, a rigorous homogenization result was then derived in \cite{ThesisPawelczyk14}.
This has a posteriori given a justification to the heuristic cell-problem computations that had been performed in 
\cite{ContiLenzRumpf2007,ContiLenzRumpf2008}  to 
study numerically the influence of the shape of the particles on the macroscopic material behavior both in composites and in polycrystals. 
Some of these results are discussed, also in view of the new homogenization result, in Section \ref{seccellpb}.
The relaxation of the model, which is appropriate for situations where the particles are much larger than the domain size, has been addressed in \cite{ContiLenzRumpf2012}.  
Our model was extended  to the study of time-dependent problems, in a setting in which each particle 
changes gradually from one phase to the other,  in \cite{ContiLenzRumpf2016}.

\section{The model}
\label{secmodel}

In this section we describe the general physical model, while the precise mathematical assumptions needed  for the homogenization result are presented in the next section.

We work in the framework of continuum mechanics and micromagnetism. For simplicity we discuss only the two-dimensional case, the extension of the model to three dimensions is immediate. 
We refer to \cite{Brown1963,HubertSchaefer} for background on micromagnetism, and to \cite{BallJames87,PitteriZanzotto02,Bhatta} for the treatment of diffusionless solid-solid phase transitions. The model we use here was first presented in 
\cite{ContiLenzRumpf2007}, to which we refer for further motivation and details.

Let $\Omega \subset \mathbb R^2$ describe the reference configuration for the composite material, and let $\omega \subset \Omega$ be the part occupied by the magnetic-shape-memory material (MSM).
We model the magnetization in the MSM-particles as a measurable function $m : \omega \to \mathbb R^2$ with $\lvert m \rvert = m_s$.
Here $m_s \in (0,\infty)$ is a parameter representing the saturation magnetization, which depends on the temperature and the specific choice of material.
For the analytical treatment we assume the temperature and material to be fixed, and thus we may assume, after normalization, $m_s = 1$.
Denoting by $S^1$ the set of unit vectors in $\mathbb R^2$, we get $m \in S^1$ on $\omega$. 
For notational convenience we extend $m$ to $\mathbb R^2 \setminus \omega$ by $0$, and introduce $S^1_0 := S^1 \cup \{ 0 \} $.

The austenite--martensite phase transition leads to the presence of $d$ ($d=2$ in our concrete application) symmetry-related variants of the martensitic phase. Magnetically, each of them has a preferred direction for the  magnetization, called easy axis; we denote those directions by
$f_1$, \dots, $f_d\in S^1$. The phase index $p(x)$ represents the variant at point $x\in\omega$. It is convenient to view it as an element of $\mathbb R^d$, taking values in a discrete set $\mathcal B:= \{e_1, e_2, \ldots, e_d\}\subset\mathbb R^d$.
This way any $\overline p$ in the convex hull of $\mathcal B$, denoted by $\mathrm{conv}\,\mathcal B$, is a unique convex combination of vectors in $\mathcal B$, allowing the tracking 
of the contribution of different variants. As before, we extend $p$ by $0$ on $\mathbb R^2 \setminus \omega$. Thus $p$ takes values in $\mathcal B_0 := \mathcal B \cup \{ 0\}\subset\mathbb R^d$.

The coupling between the phase variable $p$ and the magnetization $m$ is expressed via the anisotropy energy, which is obtained by integrating
a density $\varphi$ {depending} on $x, m(x)$, and $p(x)$ (in practice, if $x\in\omega$ then $\varphi$ 
can be seen as a function of $m(x)\cdot f_{p(x)}$). It is minimized if $x \not\in \omega$ or $m(x)$, $p(x)$ are compatible with each other, in the sense that $m(x)$ and $f_{p(x)}$ are parallel, see \cite{HubertSchaefer}.
The global anisotropy energy $E_{\mathrm{aniso}}: L^\infty( \Omega, S^1_0) \times L^\infty(\Omega, \mathcal B_0) \to [0, \infty)$ then reads
\begin{equation}\begin{aligned}\label{eq:anisotropeenergie}
E_{\mathrm{aniso}}[m,p] := \int_{\mathbb R^2} \varphi(x, m(x), p(x))\, \mathrm{d}x.
\end{aligned}\end{equation}
 A specific expression for $\varphi$ is, at this level of modeling, not needed.

The magnetization $m$  in turn induces a magnetic field, called the demagnetization field.
By Maxwell's equation this field is given by $\nabla u$, where $u$ solves 
\begin{equation}\begin{aligned}\label{eq:maxwelleq}
  \Delta u + \mathrm{div}\, m = 0 \quad \text{ on } \mathbb R^2.
\end{aligned}\end{equation}
 {Here, (\ref{eq:maxwelleq}) is} understood distributionally. 
Since we extended $m$ by zero outside $\omega$, the boundary terms are automatically included in the distributional derivatives.
The demagnetization energy $E_{\mathrm{demag}}: L^\infty(\Omega, S^1_0) \to [0,\infty)$ can be computed by
\begin{equation}\begin{aligned}\label{eq:demagenergy}
E_{\mathrm{demag}}[m] := \frac 1 2 \int_{\mathbb R^2} \lvert \nabla u_m(x)\rvert^2\, \mathrm{d}x,
  \end{aligned}\end{equation}
where $u_m$ solves $\Delta u_m + \mathrm{div}\, m = 0$~\cite[subsection 3.2.5]{HubertSchaefer} with suitable boundary conditions at infinity. 
The well-posedness of this problem will be discussed in the next section.

Next, we model the external magnetic field applied to the composite material as in~\cite[subsection 3.2.4]{HubertSchaefer}. If the external field is given by $h_e \in L^1(\Omega, \mathbb R^2)$,
then the external energy $E_{\mathrm{ext}} : L^\infty( \Omega, S^1_0) \to \mathbb R$ is 
\begin{equation}\begin{aligned}\label{eq:externalenergy}
E_{\mathrm{ext}}[m] := -\int_\Omega h_e(x) \cdot m(x)\, \mathrm{d}x.
\end{aligned}\end{equation}
This energy is {minimal for} $m$  oriented in the direction of $h_e$. 

In ferromagnets the microscopic magnetization $m$ tends to be aligned in neighbouring cells, leading to a macroscopic magnetization and to 
an energy cost for fluctuations in the magnetization itself. This behaviour if often modelled by the exchange energy $E_{\mathrm{ex}}: L^\infty(\Omega, S^1_0) \to [0,\infty]$, given by
\[
  E_{\mathrm{ex}}[m] := \varepsilon^2\int_\omega \lvert \nabla m(x) \rvert^2 \mathrm dx
\]
{for}  $m \in W^{1,2}(\omega, S^1)$ and $\infty$ otherwise (see~\cite[section 3.2.2]{HubertSchaefer}).
The parameter $\varepsilon$ measures the exchange length, which can be understood as the (small) 
length scale over which the exchange term between overlapping atomic orbitals is significant.

This concludes the description of the different components of the magnetic energy $\widetilde\Emag: L^\infty(\Omega,\mathbb R^2) \times L^\infty(\Omega,\mathcal B_0) \to [0,\infty]$, given by
\[
  \widetilde\Emag[m, p] := E_{\mathrm{aniso}}[m, p] + E_{\mathrm{demag}}[m] + E_{\mathrm{ext}}[m] + E_{\mathrm{ex}}[m].
\]
If $\omega$ consists only of sufficiently `small' grains, as we will assume in the homogenization process, then it is reasonable to assume that every grain is a single crystal, and consists only of a single domain, in which $m$ is constant. 
Indeed, in this case the exchange energy is negligible, as is shown in~\cite{De93}. Thus the final magnetic energy $\Emag: L^\infty(\Omega,\mathbb R^2) \times L^\infty(\Omega,\mathcal B_0) \to [0,\infty)$ is given by
\[
  \Emag[m, p] := E_{\mathrm{aniso}}[m, p] + E_{\mathrm{demag}}[m] + E_{\mathrm{ext}}[m].
\]

We finally address the elastic properties of the composite material. Let $\widetilde v \in W^{1,2}(\Omega, \mathbb R^2)$ be the deformation of 
the composite material, which describes the spontaneous stretching  of the MSM material in response to the external field and the subsequent deformation of the polymer.
 Let $\widetilde W = W(x, F, p)$ be the (nonlinear) elastic energy density,  where $F(x)=\nabla v(x)$ is the deformation gradient at $x$, and $p(x)\in\{0,1,\dots, d\}$ denotes the martensitic variant (as usual, $p=0$ in the polymer).
 A possible choice would be $\widetilde W(x, F, p) = \mathrm{dist}^2(F R(x), \mathrm{SO}(2)(\mathrm{Id} - A_p))$ if $x \in \omega$ and $\widetilde W(x, F, p) = \mathrm{dist}^2(F, \mathrm{SO}(2))$ else.
Here the rotation $R\in L^\infty(\omega,\mathrm{SO}(2))$ encodes the local orientation of the crystal lattice, and $A_p$  the eigenstrain of the phase $p$.
The magnetization in the deformed configuration,  where Maxwell's equation has to hold, is then  $m \circ \widetilde v$.
Thus, we now would need to solve $\Delta u + \mathrm{div}\, ( m \circ\widetilde  v ) = 0 $ on $\mathbb R^2$. In addition, the anisotropy energy $\varphi$ {has} to take into account the fact that the local magnetization $m(x)$ needs to be compared with the easy axis in the deformed (and, possibly, rotated) crystal.

However, experiments show that the displacements are moderate, and thus we approximate $\widetilde v(x) = x + v(x)$
with a small displacement $v$. Consistently expanding all terms to the first nontrivial order in $m$ and $v$, we consider the magnetization $m$ instead of $m\circ\widetilde v$ and ignore elastic rotations in the anisotropy energy.
Therefore we will work with the energy
\[
  \Evoll[v,m,p] := E_{\mathrm{def}}[v,m,p] + \Emag[m,p],
\]
where $E_{\mathrm{def}}[v,m,p]$ is the linear elastic energy stored in the displacement $v$ {and} obtained integrating a quadratic form $W$  {corresponding} to the linearization of $\widetilde W$.

\section{Analytical homogenization}
\label{sechomo}

The MSM--polymer composite is a finely microstructured material, and therefore very difficult to simulate directly. A direct discretization would need to resolve explicitly
all scales from the one of the MSM particles, of the order of a few microns, to the one of the macroscopic sample, which may be in the centimeter range. The theory of homogenization
permits both a simpler qualitative understanding of the material behavior and an efficient simulation. The key idea is to exploit the difference in length scales to address the two scales separately.
Under specific assumptions of the microstructure (periodicity in the present case) one derives an effective macroscopic material model, which describes the behavior of the composite on a length scale much larger than
the one of the microstructure, by solving suitable
cell problems. At this level one does not need to address the macroscopic shape of the sample, boundary conditions or applied forces, but instead studies only the microstructure.
One then uses, in a separate process, the effective material model to perform macroscopic simulations, {which takes into account}  boundary conditions and external forces.
Typical examples of homogenization include the process leading from atomistic models of matter to continuum elasticity and the treatment of materials obtained by mixing two different components on a very fine scale.
Whereas the general theory for composites of two materials with different elastic properties is well developed, we are here interested in a situation in which the magnetic properties are also different, and coupled to an elastic phase transition.
A purely magnetic problem was studied in \cite{Pisante2004}.
The homogenization result for our system does not follow from the general theory, but needs to be proven specifically for the case at hand.
For a comprehensive overview over the subject, we refer to~\cite{Braides98,DonatoCioranescu,AllaireBuch,Milton2002}.

\begin{figure}
 \begin{center}
 \includegraphics[width=5cm]{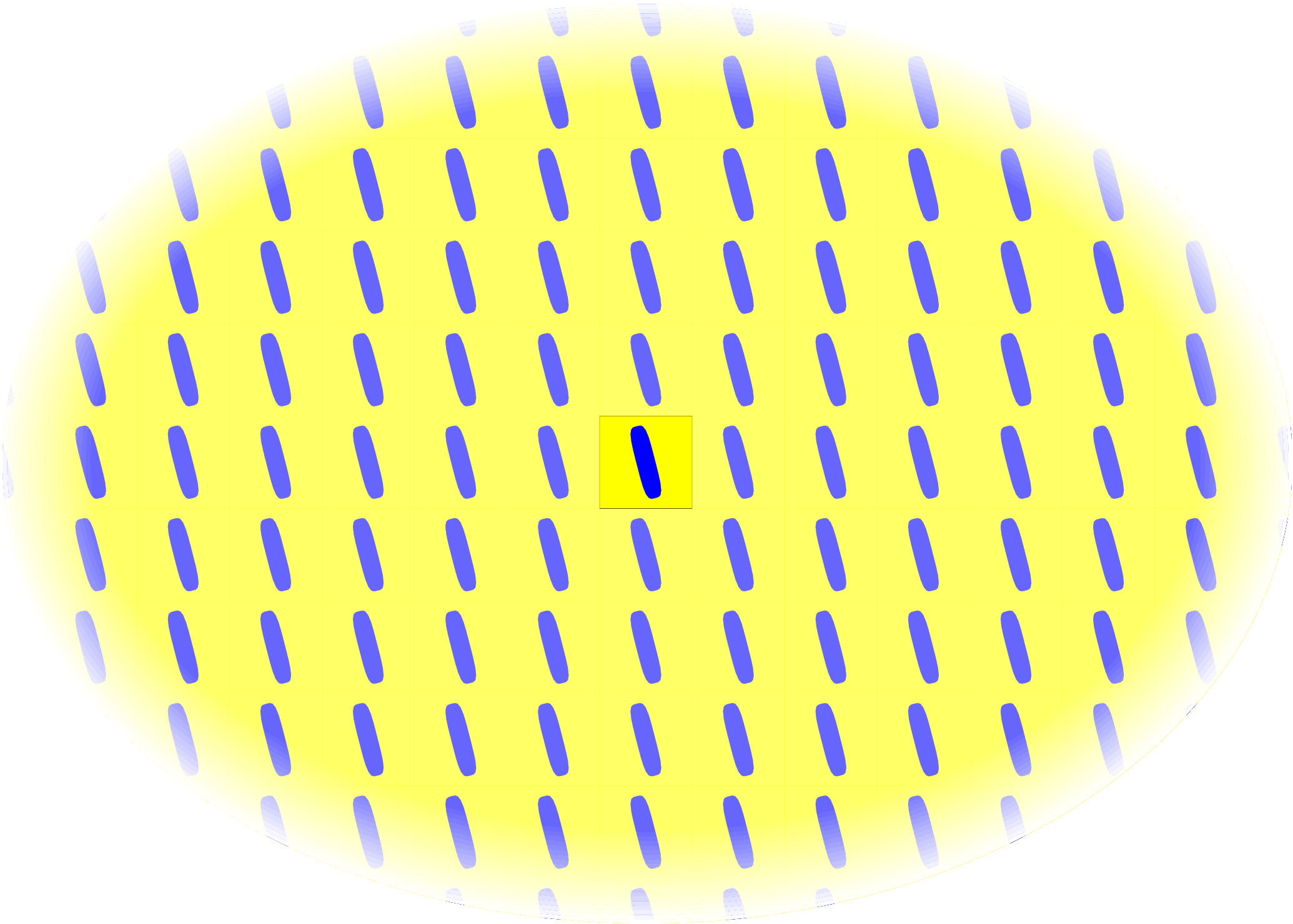}\hskip1cm
 \includegraphics[width=5cm]{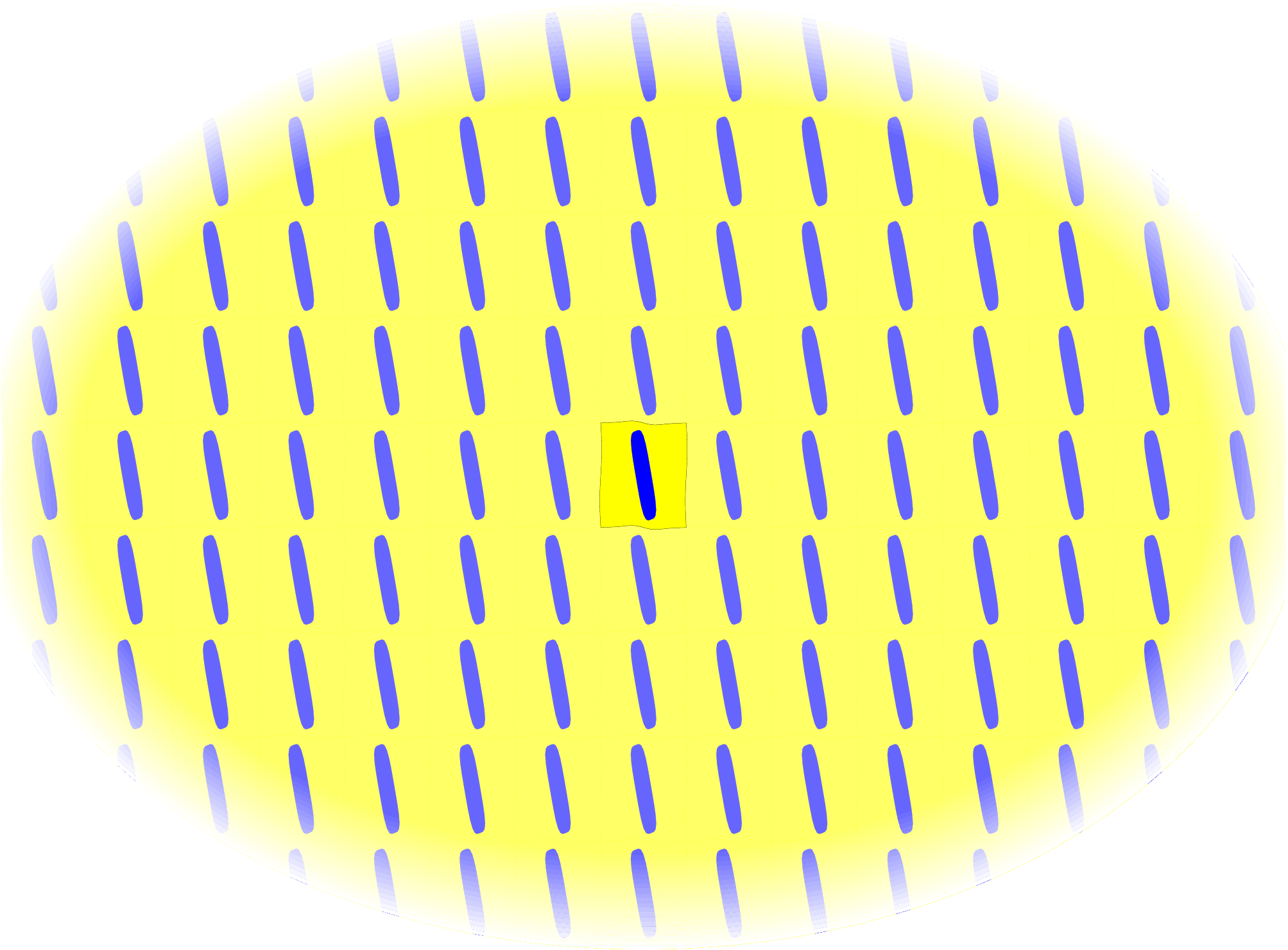}
  \end{center}
\caption{Sketch of the two-scale approach to the behavior of composites. The periodic cell is emphasized in the center.
Left picture: reference configuration. Right picture: deformed configuration. Both the particle and the shape of the unit cell are deformed.}
\label{fighomosketch}
\end{figure}

Let $\Omega \subset \mathbb R^2$ be an open, bounded Lipschitz domain, and let $\omega \subset Q:=(0,1)^2$ be measurable.
For every $k\in \mathbb N$ let $\omega_k := \Omega \cap \big( \frac 1 k (\omega + \mathbb Z^2) \big )$, see Figure~\ref{fighomosketch} for an illustration.
It is easily seen that $\lim_{k\to \infty}\mathcal L^2( \omega_k) = \mathcal L^2(\omega) \,\mathcal L^2(\Omega)$  is the volume fraction of the composite  occupied by the MSM particles.
Here, $\mathcal L^2$ is the Lebesgue measure.
The set of admissible magnetization fields and phase indices for $\omega_k$ is given by
\[
  S^{\mathrm{mag}}_k := \left \{ (m,p) \in L^\infty(\mathbb R^2, \mathbb R^2) \times L^\infty(\mathbb R^2, \mathbb R^d): \; \lvert m \rvert = \lvert p\rvert = \chi_{\omega_k} \right \},
\]
where $\chi_{\omega_k}$ is the characteristic function of $\omega_k$.

In order to compute the demagnetization energy we need to show the well-posedness of~\eqref{eq:maxwelleq}.
By Lax-Milgram's Theorem there exists for any measurable $m$ with $\lvert m \rvert = \chi_{\omega_k}$ a unique solution $u$ of~\eqref{eq:maxwelleq} in the space
\begin{align*}
  W^{1,2}_*(\mathbb R^2) := \left \{ u \in W^{1,2}_{\mathrm{loc}}(\mathbb R^2): \int_{\{\lvert x\rvert \leq 1\}}\hspace{-0.3cm} u(x) \mathrm dx = 0 \; \text{ and } \; \nabla u \in L^2(\mathbb R^2, \mathbb R^2) \right \}.
\end{align*}
The condition on the average of $u$ on the unit ball is an arbitrarily chosen criterion used to make the solution unique, since both the equation defining $u$ and the energy only depend on {the} gradient {of $u$}. Alternatively one could work in the space of curl-free $L^2$ vector fields. The explicit appearance of the potential $u$ in the definition of the space is useful for truncation and interpolation.
We also remark that, due to the fact that the fundamental solution of Laplace's equation is logarithmic in two dimensions, the potential $u$ is not in $L^2(\mathbb R^2)$. In  $\mathbb R^3$ instead one would obtain a unique solution $u\in W^{1,2}(\mathbb R^3)$.
Given a magnetization field $m$, we denote by  $u_m$ the solution of~\eqref{eq:maxwelleq} in $W^{1,2}_*(\mathbb R^2)$.
Let the anisotropic energy density $\varphi: \mathbb R^2 \times S^1_0 \times \mathcal B_0 \to [0,\infty)$ in~\eqref{eq:anisotropeenergie}
satisfy
\begin{enumerate}
\item $\varphi$ is a Carath\'eodory function, i.e.,
\begin{enumerate}
  \item $\varphi(\cdot, m, p)$ is measurable for all $(m,p) \in S^1_0 \times \mathcal B_0$. 
  \item  $\varphi(y, \cdot, \cdot)$ is continuous for all $y \in \mathbb R^2$.
\end{enumerate}
\item The map $y \mapsto \varphi(y,m, p)$ is $Q$-periodic for any $m\in S^1_0, p\in\mathcal B_0$.
\item For any $y\in \mathbb R^2 \setminus \omega$ we have $\varphi(y, \cdot, \cdot) = 0$.
\item There exists $M > 0$, such that for all $m, m' \in S^1, \, p \in\mathcal B, \; y\in \omega$ we have
\[
  \lvert \varphi(y, m, p ) - \varphi(y, m', p) \rvert \leq M \lvert m - m' \rvert.
\]
\end{enumerate}
For the external magnetic field we  assume $h_e \in L^1(\Omega, \mathbb R^2)$. 
We define the magnetic energy $\Emag_k: L^\infty(\mathbb R^2, \mathbb R^2) \times L^\infty(\mathbb R^2, \mathbb R^d)\to (-\infty, \infty]$ by
\begin{align*} 
  \Emag_k[m,p] := \int_\Omega \varphi(kx, m(x), p(x))\mathrm dx + \frac 1 2 \int_{\mathbb R^2} \lvert \nabla u_m \rvert^2 \mathrm dx - \int_\Omega m(x) \cdot h_e(x) \mathrm d x
\end{align*}
if $(m,p) \in S_k$,
and $\Emag_k[m,p] =\infty$ otherwise.

The homogenization limit corresponds to taking weak limits of $p_k$ and $m_k$, which are $L^\infty$ fields.
It is therefore not surprising that the conditions $|p_k|=|m_k|=\chi_{\omega_k}$ will be relaxed in the limit, with $\chi_{\omega_k}$ converging weakly to the 
function  $\mathcal L^2(\omega)\chi_\Omega$, and the equality conditions being replaced by  inequalities representing the convexification of the conditions. Indeed, we will show that
the space of the possible limiting values $m$ and $p$  is
\begin{align*} 
S^{\mathrm{mag}}_\infty :=&  \{ (m,p) \in L^\infty(\mathbb R^2, \mathbb R^2) \times L^\infty(\mathbb R^2, \mathrm{conv}\,\mathcal B_0):\\
  &\hspace{1cm} \lvert m \rvert \leq \mathcal L^2(\omega) \chi_{\Omega}, \; \lvert p\rvert_{\ell^1} = \mathcal L^2(\omega) \chi_{\Omega}  \},
  \end{align*}
where $\lvert p\rvert_{\ell^1} = \sum_{k=1}^d \lvert p_k\rvert$ and $\mathrm{conv}\,\mathcal B_0$ is the convex hull of $\mathcal B_0$.

The effective energy $\Emag_\infty: L^\infty(\mathbb R^2, \mathbb R^2) \times L^\infty(\mathbb R^2, \mathbb R^d)\to (-\infty, \infty]$ on $S^{\mathrm{mag}}_\infty$ is then given by 
\begin{align*} 
  \Emag_\infty[m,p] := \int_\Omega f^{\mathrm{mag}}_{\mathrm{hom}}(m(x), p(x))\mathrm dx + \frac 1 2 \int_{\mathbb R^2} \lvert \nabla u_m \rvert^2 \mathrm dx - \int_\Omega m(x) \cdot h_e(x) \mathrm d x,
\end{align*}
and $\Emag_\infty = \infty$ otherwise. 
Here $f^{\mathrm{mag}}_{\mathrm{hom}} : \mathbb R^2 \times \mathrm{conv}\, \mathcal B_0 \to [0, \infty]$ is given by
\begin{equation}\begin{aligned}\label{eq:fhomdefinition}
f^{\mathrm{mag}}_{\mathrm{hom}}(\mu,\rho) = \inf_{ ((M_k, P_k))_{k \in \mathbb N} \in \mathcal A_{\mu, \rho} } \liminf_{k \to \infty} \left(\int_Q \varphi(kx, M_k(x), P_k(x))+ \frac 1 2\lvert \nabla U_{M_k} \rvert^2 \mathrm dx \right)
\end{aligned}\end{equation}
where $U_M \in W^{1,2}_{\mathrm{loc}}(\mathbb R^2, \mathbb R^2)$ is the unique function satisfying
\[
  U_M \text{ is $Q$-periodic}, \; \int_Q U_M(x)\mathrm dx = 0 \;  \text{ and } \; \Delta U_M + \mathrm{div}\, M = 0 \; \text{ on } Q,
\]
and
\begin{align*} 
  \mathcal A_{\mu, \rho} := & \Big\{ ((M_k , P_k))_{k \in \mathbb N} \subset L^\infty(Q, S^1_0) \times L^\infty(Q,  \mathcal B_0): \\
                            &  \quad \forall k \in \mathbb N: \lvert M_k \rvert = \lvert P_k \rvert = \chi_{\omega}, \quad \int_Q M_k \mathrm dx= \mu \quad \text{ and } \\\
                            & \quad M_k \stackrel*\rightharpoonup \mu \text{ on } L^\infty(Q, \mathbb R^2), \quad P_k \stackrel*\rightharpoonup \rho \text{ on } L^\infty(Q, \mathbb R^d)\Big\}.
\end{align*}

We now can state the homogenization {theorem}  for the purely magnetic energy.
This and the next results have been proven in \cite{ThesisPawelczyk14} extending the classical homogenization theory for elastic materials \cite{Braides98,DonatoCioranescu,AllaireBuch,Milton2002}
and the homogenization results for three-dimensional purely magnetic materials in \cite{Pisante2004}; a brief summary of the proof is given below.
\begin{thm}
With the above assumptions on $\varphi$ we have the following results:
\begin{enumerate}
\item\label{enum:compactness} Let $((m_k, p_k))_{k \in \mathbb N} \subset L^\infty( \mathbb R^2, \mathbb R^2) \times L^\infty(\mathbb R^2, \mathbb R^d)$ with
$$\liminf_{k\to \infty} \Emag_k[m_k,p_k] < \infty\,.$$
Then there exist a subsequence (not relabeled), and $(m,p) \in S^{\mathrm{mag}}_\infty$ such that 
\begin{align*} 
&m_k \stackrel*\rightharpoonup m \text{ in } L^\infty(\Omega,\mathbb R^2), 
\quad p_k \stackrel*\rightharpoonup p \text{ in } L^\infty(\Omega,\mathbb R^d), \\
&\nabla u_{m_k} \rightharpoonup \nabla u_m \text{ in } L^2(\mathbb R^2),
\quad \nabla u_{m_k} \to \nabla u_m \text{ in } L^2(\mathbb R^2\setminus \Omega),
  \end{align*}
where $u_m, u_{m_k}\in W^{1,2}_*(\mathbb R^2)$ are solutions to 
\[
  \Delta u_{m_k} + \mathrm {div}\, {m_k} = \Delta u_m + \mathrm {div}\, m = 0 \;\text{on $\mathbb R^2$.}
\] 
\item\label{enum:Gammalimsup} For any $(m,p)\in L^\infty(\mathbb R^2, \mathbb R^2) \times L^\infty(\mathbb R^2, \mathbb R^d)$ \\ there exists a sequence $((m_k,p_k))\subset L^\infty(\mathbb R^2, \mathbb R^2) \times L^\infty(\mathbb R^2, \mathbb R^d)$
  such that 
  \[
    \limsup_{k\to \infty} \Emag_k[m_k, p_k] = \Emag_\infty[m,p].
  \]
\item\label{enum:Gammaliminf} Let $((m_k,p_k))\subset L^\infty(\mathbb R^2, \mathbb R^2) \times L^\infty(\mathbb R^2, \mathbb R^d)$ \\ with $m_k \stackrel*\rightharpoonup m, p_k \stackrel*\rightharpoonup p$ in $L^\infty(\Omega, \mathbb R^2), L^\infty(\Omega,\mathbb R^d)$ resp. Then
\[
    \liminf_{k\to \infty} \Emag_k[m_k, p_k] \geq \Emag_\infty[m,p].
\]
\end{enumerate}
\end{thm}
As usual, this result can be easily stated in terms of $\Gamma$-convergence.

The compactness result~\ref{enum:compactness} for $m_k, p_k$ follows easily, since the uniform bound on the energy implies that there is a subsequence with  $\lvert m_k \rvert = \lvert p_k \rvert = \chi_{\omega_k}$ almost everywhere. 
Thus the $L^\infty$-norm of $m_k, p_k$ is uniformly bounded by $1$, and we can extract a converging subsequence. The convergence of the demagnetization field then follows from~\cite[Lemma~A.5.4]{ThesisPawelczyk14}.

A key ingredient in the proof is the fact that 
 $f^{\mathrm{mag}}_{\mathrm{hom}}$ can be equivalently defined by minimizing over other sets than $\mathcal A_{\mu,\rho}$.
Specifically, we can either drop  the condition of having a fixed mean value on $m_k$ on $Q$, or drop the condition that $M_k \stackrel*\rightharpoonup \mu$ without changing the value of $f^{\mathrm{mag}}_{\mathrm{hom}}(\mu, \rho)$.
Also instead of solving $\Delta U + \mathrm{div}\, M = 0$ on $Q$ with periodic boundary, we can solve it on $Q$ with vanishing Dirichlet boundary condition. Furthermore the density $f^{\mathrm{mag}}_{\mathrm{hom}}$ does not change if  $U$ 
is replaced by the solution 
$U \in W^{1,2}_*(\mathbb R^2)$ with $\Delta U + \mathrm{div}\, \chi_{Q}M = 0$ and the term $\int_Q \lvert \nabla U \rvert^2$ in~\eqref{eq:fhomdefinition} replaced by the integral $\int_{\mathbb R^2} \lvert \nabla U\rvert^2$.
We refer to \cite[Proposition~2.2.1]{ThesisPawelczyk14} for a precise statement of these facts and a proof.

This equivalence is crucial to relate the macroscopic demagnetization field $u$ to the microscopic fields $U$, which are the result of oscillations in the magnetic charges.
More precisely, we first prove that $f^{\mathrm{mag}}_{\mathrm{hom}}$ is a continuous function and $\Emag_\infty$ restricted to the limiting function $S^{\mathrm{mag}}_\infty$ is continuous w.r.t.\ the $L^1$-topology.
Thus it suffices to prove~\ref{enum:Gammalimsup} for step functions $m,p$. This is done by covering $\Omega$ by squares, and assuming $m,p$ are constant on each of them. On each square we can explicitly construct a recovery sequence {of} periodic functions, and by glueing them together we obtain a global recovery sequence.
Details are given in~\cite[Theorem 2.3.1 and Theorem 2.3.2]{ThesisPawelczyk14}.

\vspace{\baselineskip}
We now include the elastic energy in the energy in terms of a displacement $v$. 
Let $W : \mathbb R^2 \times \mathbb R^{2\times 2} \times \mathcal B_0 \to \mathbb R$ be a Carath\'eodory function, such that
\begin{enumerate}
  \item $W$ is $Q$-periodic in the first component. 
  \item There exists $r: [0,\infty) \to [0, \infty)$ with $r(\delta) \searrow 0$ for $\delta \searrow 0$, 
        such that for any $x \in \mathbb R^2, \, F,G \in\mathbb R^{2\times 2}, \, p \in \mathcal B_0$ it holds
        \[
            W(x, F +G, p) \leq ( 1 + r(\lvert G\rvert)) W(x,F,p) + r(\lvert G\rvert).
        \]
  \item There exists $c > 0$, such that for any $x \in \mathbb R^2, \, F,G \in\mathbb R^{2\times 2}, \, \rho \in \mathcal B_0$ it holds
        \[
            \frac { \lvert \mathrm{sym}\, F \rvert^2} c - c\leq W(x, F, p) \leq c\lvert \mathrm{sym}\, F \rvert^2  + c,
        \]
      {where $\mathrm{sym}\, F =\tfrac12 (F+F^T)$.}
\end{enumerate}
We define the set of admissible triples $(v,m,p)$ by
\[
  S_k := \left \{ (v, m,p) \in W^{1,2}(\Omega, \mathbb R^2) \times S_k^{\mathrm{mag}} \right\},
\]
and the energy $\Evoll_k: W^{1,2}(\Omega,\mathbb R^2) \times  L^\infty(\mathbb R^2, \mathbb R^2) \times L^\infty(\mathbb R^2, \mathbb R^d)\to (-\infty, \infty]$ by
\begin{align*}
  \Evoll_k[v, m,p] :=   \Emag_k[m,p] + \int_\Omega W(kx, \nabla v(x), p(x))\mathrm dx
\end{align*}
 if $(v,m,p) \in S_k$,
and $\Evoll_k[v, m,p] =\infty$ otherwise.
The set of limiting configurations is given by $S_\infty := W^{1,2}(\Omega, \mathbb R^2) \times S_\infty^{\mathrm{mag}}$, and the effective energy is 
\[
  \Evoll_\infty[v,m,p] := \int_\Omega f_{\mathrm{hom}}(\nabla v(x), m(x), p(x)) \mathrm dx + \frac 1 2 \int_{\mathbb R^2} \lvert \nabla u_m \rvert^2 \mathrm dx - \int_\Omega m \cdot h_e \mathrm dx,
\]
where 
\begin{align*}
&f_{\mathrm{hom}} (\beta, \mu, \rho):= \inf_{ ((V_k, M_k, P_k)) \in \mathcal A_{\beta, \mu,\rho} } \liminf_{k\to\infty}  \\
&\quad\left(\int_Q W(kx, \nabla V_k(x),P(x))+ \varphi(kx, M_k(x), P_k(x))+ \frac 1 2\lvert \nabla U_{M_k} \rvert^2 \mathrm dx \right),
\end{align*}
and
\begin{align*}
  \mathcal A_{\beta,\mu, \rho} := & \Big\{ ((V_k, M_k , P_k))_{k \in \mathbb N} \subset W^{1,2}(Q, \mathbb R^2) \times L^\infty(Q, S^1_0) \times L^\infty(Q,  \mathcal B_0): \\
                                    & \quad ((M_k, P_k)) \in \mathcal A_{\mu,\rho}, \;\; \nabla V_k \rightharpoonup \beta \text{ in } L^2(Q,\mathbb R^2)  \Big\}.
\end{align*}

\begin{thm}
With the above assumptions on $\varphi$ and $W$ we have the following results:
\begin{enumerate}
\item\label{enum:compactnessfull} Let $((v_k, m_k, p_k))_{k \in \mathbb N} \subset W^{1,2}(\Omega, \mathbb R^2)\times L^\infty( \mathbb R^2, \mathbb R^2) \times L^\infty(\mathbb R^2, \mathbb R^d)$ \\ with $\liminf_{k\to \infty} \Evoll_k[m_k,p_k] < \infty$. \\
Then there exist a subsequence (not relabeled), and $(v, m,p) \in S_\infty$ such that 
\begin{align*} 
&m_k \stackrel*\rightharpoonup m \text{ in } L^\infty(\Omega,\mathbb R^2), 
\quad p_k \stackrel*\rightharpoonup p \text{ in } L^\infty(\Omega,\mathbb R^d), \\
&\nabla u_{m_k} \rightharpoonup \nabla u_m \text{ in } L^2(\mathbb R^2),
\quad \nabla u_{m_k} \to \nabla u_m \text{ in } L^2(\mathbb R^2\setminus \Omega), \\
&\left(v_k - \frac 1 {\mathcal L^2(\Omega)} \int_\Omega v_k\right) \rightharpoonup v \text{ in } W^{1,2}(\Omega, \mathbb R^2),
  \end{align*}
where $u_m, u_{m_k}\in W^{1,2}_*(\mathbb R^2)$ are solutions of
\[
  \Delta u_{m_k} + \mathrm {div}\, {m_k} = \Delta u_m + \mathrm {div}\, m = 0 \;\text{on $\mathbb R^2$.}
\] 
\item\label{enum:Gammalimsupfull} For any $(v, m,p)\in W^{1,2}(\Omega,\mathbb R^2) \times L^\infty(\mathbb R^2, \mathbb R^2) \times L^\infty(\mathbb R^2, \mathbb R^d)$ there exists a sequence $((v_k,m_k,p_k))\subset W^{1,2}(\Omega,\mathbb R^2)\times L^\infty(\mathbb R^2, \mathbb R^2) \times L^\infty(\mathbb R^2, \mathbb R^d)$
  such that 
  \[
    \limsup_{k\to \infty} \Evoll_k[v_k, m_k, p_k] = \Evoll_\infty[v, m,p].
  \]
\item\label{enum:Gammaliminffull} Let $((v_k,m_k,p_k))\subset W^{1,2}(\Omega,\mathbb R^2) \times L^\infty(\mathbb R^2, \mathbb R^2) \times L^\infty(\mathbb R^2, \mathbb R^d)$ with $m_k \stackrel*\rightharpoonup m, \,p_k \stackrel*\rightharpoonup p$ in $L^\infty(\Omega, \mathbb R^2), L^\infty(\Omega,\mathbb R^d)$ resp., and $v_k \rightharpoonup v$ in $W^{1,2}(\Omega, \mathbb R^2)$. Then
\[
    \liminf_{k\to \infty} \Evoll_k[v_k, m_k, p_k] \geq \Evoll_\infty[v,m,p].
\]
\end{enumerate}
\end{thm}
The proof is similar to the previous one. For~\ref{enum:compactnessfull} it only remains to prove compactness for the sequence $(v_k)$. 
The coercivity of $W$ yields an uniform $L^2$-bound on the sequence $(\mathrm{sym}\, \nabla v_k)$. An application of Korn's and Poincar\'e's inequality yields then compactness.

While the proof of~\ref{enum:Gammaliminffull} is analogous to the previous one, we need to be more careful in~\ref{enum:Gammalimsupfull}.
To approximate $v$ in $W^{1,2}$ we use piecewise linear functions on triangles. This is, however, not the natural decomposition for the $Q$-periodic microstructure. This can be dealt with by covering each triangle by squares once more.
Detailed proofs can be found in~\cite[Theorem 3.3.1 and Theorem 3.3.2]{ThesisPawelczyk14}.

\section{Numerical study of the cell problem}
\label{seccellpb}

We illustrate the significance of the homogenization result by studying numerically the cell problem for a few selected microstructures. Our results characterize
the dependence of the macroscopic material properties on the microstructure. From the viewpoint of applications, the key interest lies in devising microstructures
which can be practically produced and which lead to {(near)} optimal macroscopic properties, such as for example a large transformation strain in response to an external magnetic field or a large work
output. This constitutes a shape-optimization problem:
 we are interested in devising the shape and location of the particles which optimizes some scalar quantity. It is however  different from most classical 
shape-optimization problems, in that it completely takes place at the microstructural level, and in that only a very restricted set of shapes is practically relevant. Indeed, the control on the shape and the ordering
of the MSM particles is only very indirect, as they are obtained for example by mechanically grinding thin ribbons of MSM material \cite{LiuScheerbaumWeissGutfleisch2009}.  Their orientation can, up to a certain point, be controlled by applying a magnetic field during the solidification of the polymer.  At the same time, the elastic properties of the polymer can be, up to a  {certain degree}, tuned by choosing different compositions.
We therefore select a few geometric and material parameters (particle elongation, Young modulus of the polymer, etc.)  and investigate their influence on the macroscopic material behavior. Our results can serve as a guide for the experimental search for production techniques which lead to composites with {successively improved}  properties.

As customary in the theory of homogenization, we assume that the corrector fields  $(V_k, M_k, P_k)$ 
entering $f_{\mathrm{hom}} (\beta, \mu, \rho)$
have the same periodicity of the microstructure (that is, $(0,1/k)^2$) and that $V_k$ obeys affine-periodic boundary conditions. In practice, 
we  approximate
$f_{\mathrm{hom}}(\beta, \mu, \rho)$ by $f_{\mathrm{hom}}^{(1)}(\beta, \mu, \rho)$, which (after rescaling) is defined as the infimum of 
\begin{align*}
\int_Q W(x, \nabla V(x),P(x))+ \varphi(x, M(x), P(x))+ \frac 1 2\lvert \nabla U_{M} \rvert^2 \mathrm dx
\end{align*}
over all $V\in  W^{1,2}(Q, \mathbb R^2) $ which obey $V(x+e_i)=V(x)+\beta e_i$ in the sense of traces for $x\in\partial Q$, all
$M\in L^\infty(Q, S^1_0)$ which obey $|M|=\chi_\omega$ and $\int_Q M dx = \mu$, 
and all $P\in L^\infty(Q,  \mathcal B_0)$ with $|P|=\chi_\omega$. As usual in nonconvex homogenization, the usage of larger unit cells may lead to spontaneous symmetry breaking and formation of microstructures on intermediate scales, as {it} was discussed in \cite[Sect. 6.3 and Fig. 13]{ContiLenzRumpf2007}. For the sake of simplicity we do not address this issue here.
We instead include the magnetostrictive effect by solving Maxwell's equations in the deformed configuration (an effect that is not included in the first-order model used in the homogenization).

We use a boundary-element method for solving both the linear elastic problem (with piecewise affine ansatz functions) and the magnetic problem (with piecewise constant ansatz functions). Both the particle-matrix interfaces and the boundary of the unit cell are approximated by polygons.  Each particle is assumed to have a uniform magnetization and to deform affinely; the boundary conditions on $\partial(0,1)^2$ are affine-periodic, as usual in the theory of homogenization. Our numerical algorithm is based on computing the (discrete) gradient of the energy in the relevant variables, performing a line search in the {descent} direction, and then updating the direction as in the conjugate gradient scheme.

We choose parameters that match the experimentally known values for NiMnGa particles,
as in \cite{ContiLenzRumpf2007,ContiLenzRumpf2016}.
In  the magnetic energy we use $\frac{M_s}{\mu_0} = 0.50 \frac{\rm MPa}{\rm T}$, $\frac{M_s^2}{\mu_0} = 0.31 \,\rm MPa$,
$K_u= 0.13\,\rm MPa$. We assume an external magnetic field of $1\, \rm T$.
The elastic constants {taken into account} for NiMnGa are $\epsilon_0 = 0.058$, $C_{11} = 160 \,\rm GPa$, $C_{44}= 40 \,\rm GPa$, $C_{11}-C_{12} = 4 \, \rm GPa$.
The elastic modulus of the polymer varies between $E= 0.03\, \rm MPa$ and $E= 80\, \rm MPa$, its Poisson ratio is assumed fixed at $\nu= 0.45$. 

We assume here that the polymer is solidified with the MSM particles in one of the two martensitic variants. This introduces a clear preference for one of the two, and renders transformation more difficult. At the same time it introduces a natural mechanism for transforming back to the original state after the external field is removed. From the viewpoint of material production, the current setting corresponds to the assumption that the solidification of the polymer occurs below the critical temperature for the solid-solid phase transformation (which is about $70^\circ$ C). We remark that in our previous work \cite{ContiLenzRumpf2007} we had instead assumed that solidification occurs in the austenitic phase, rendering the two martensitic variants symmetric.

\begin{figure}
 \begin{center}
  \begin{minipage}[b]{1.6cm}
   \includegraphics[width=1.9cm]{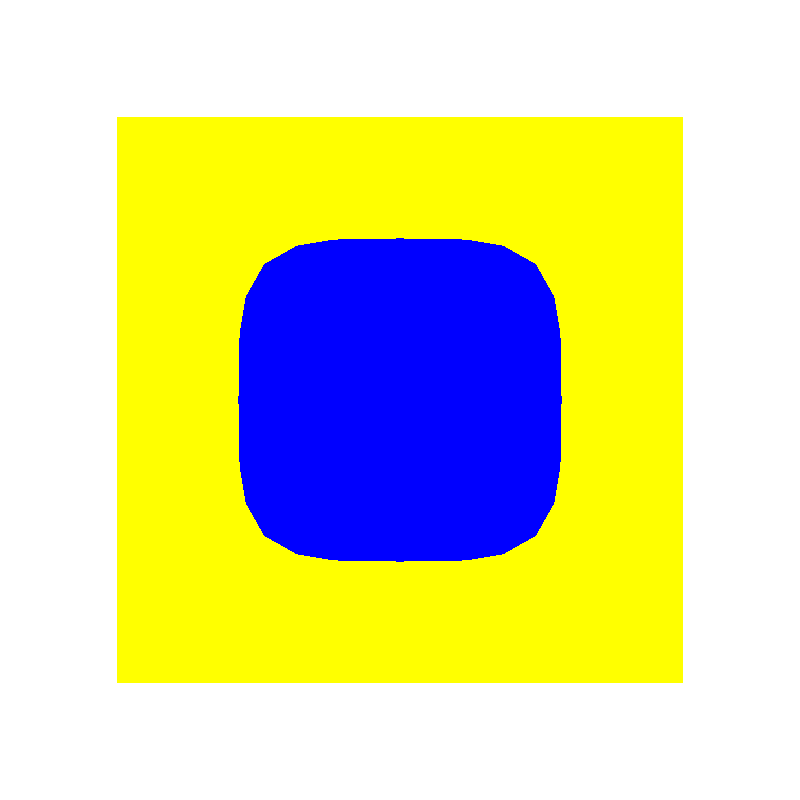}\\
   \includegraphics[width=1.9cm]{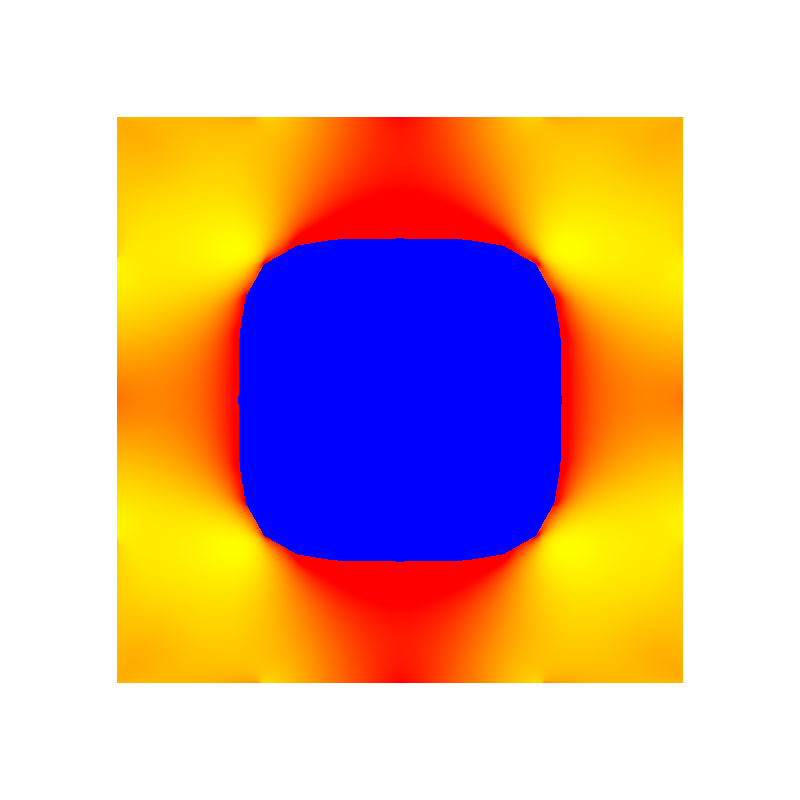}\\
  \end{minipage}\hspace{1mm}
  \includegraphics[width=4.7cm]{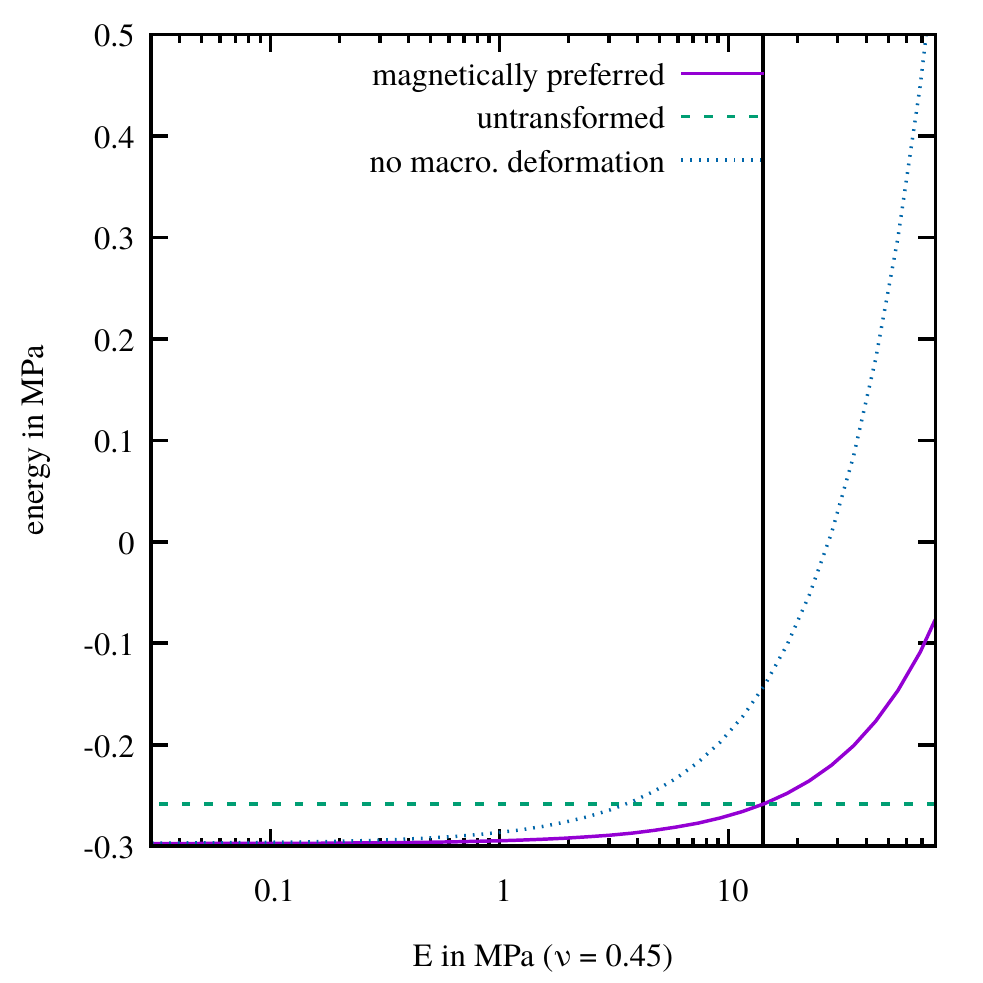}\hskip4mm
  \includegraphics[width=4.7cm]{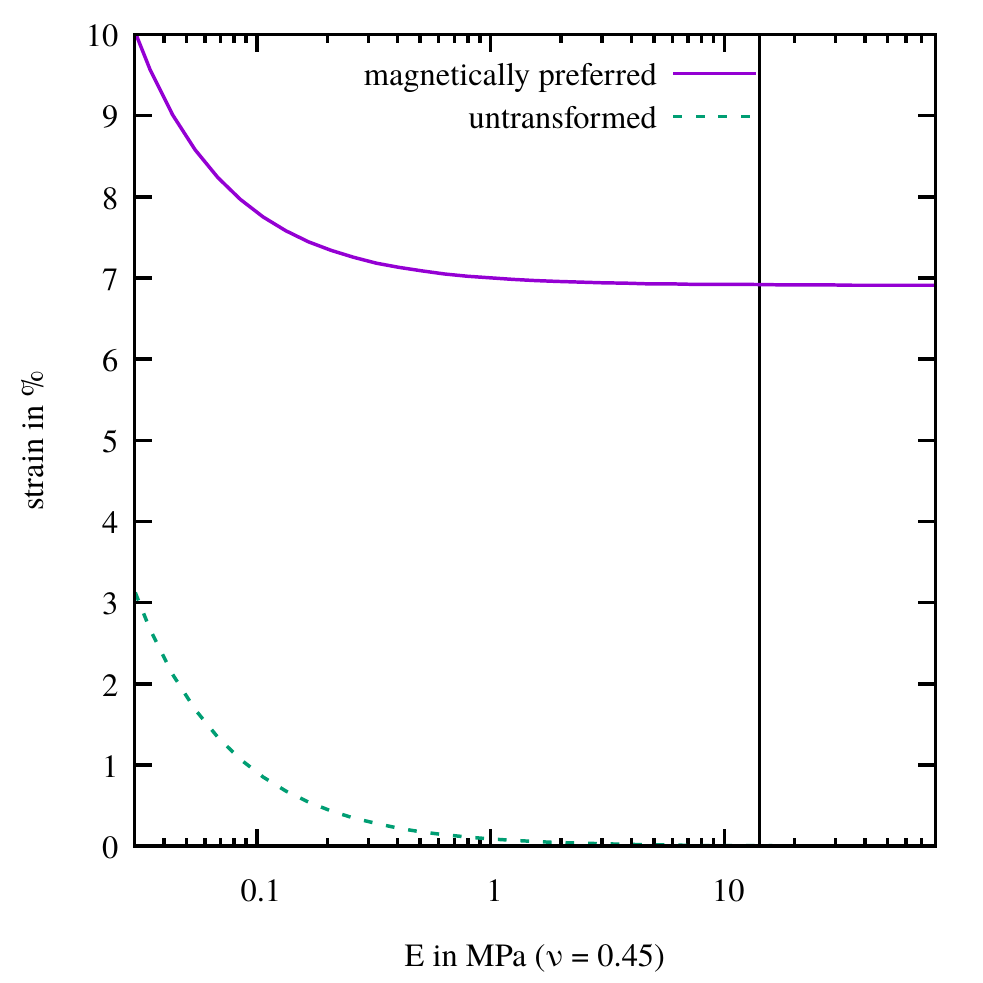}
 \end{center}
\caption{Energy and macroscopic spontaneous strain as a function of the polymer's elastic modulus, for the two martensitic variants (solid line: magnetically preferred variant; dashed line: untransformed). The left plot depicts in addition the energy if no macroscopic deformation is allowed (dotted line). The sketches on the far left show the geometric configuration of the particle and the elastic energy density in the polymer (yellow: 0 kPa, red: 0.1 kPa, for an elastic modulus of 1 MPa).}
\label{figEmodul}
\end{figure}

\begin{figure}
 \begin{center}
  \begin{minipage}[c]{1.3cm}
   \includegraphics[width=1.3cm]{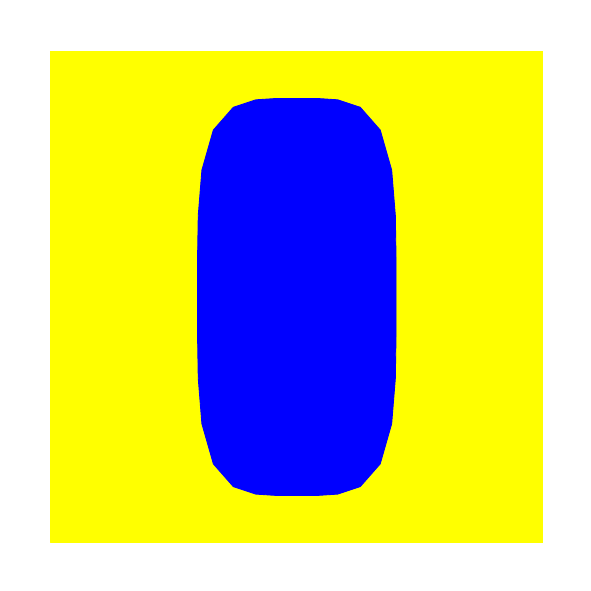}\\
   \includegraphics[width=1.3cm]{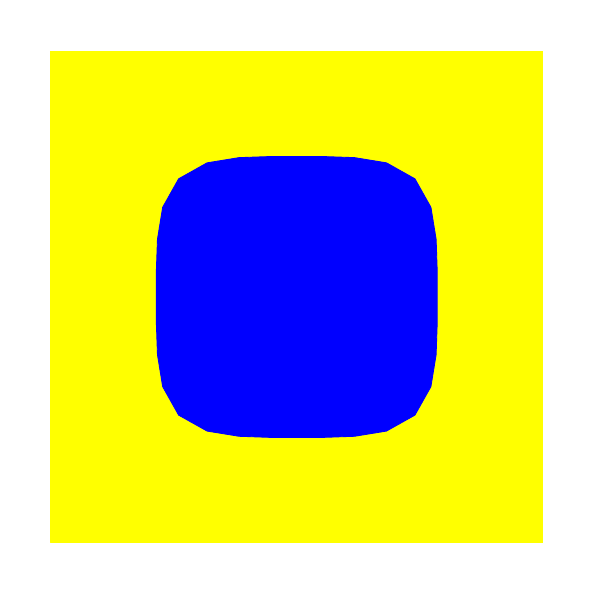}\\
   \includegraphics[width=1.3cm]{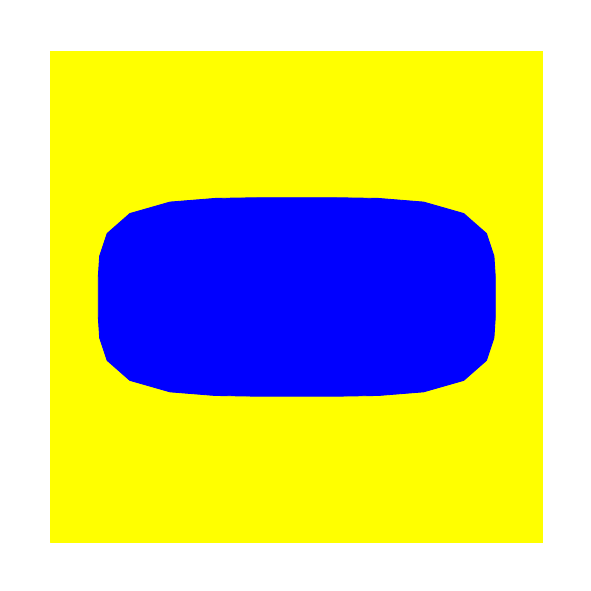}\\[0mm]
  \end{minipage}\hspace{4mm}
  \begin{minipage}[c]{9.8cm}
   \includegraphics[width=4.7cm]{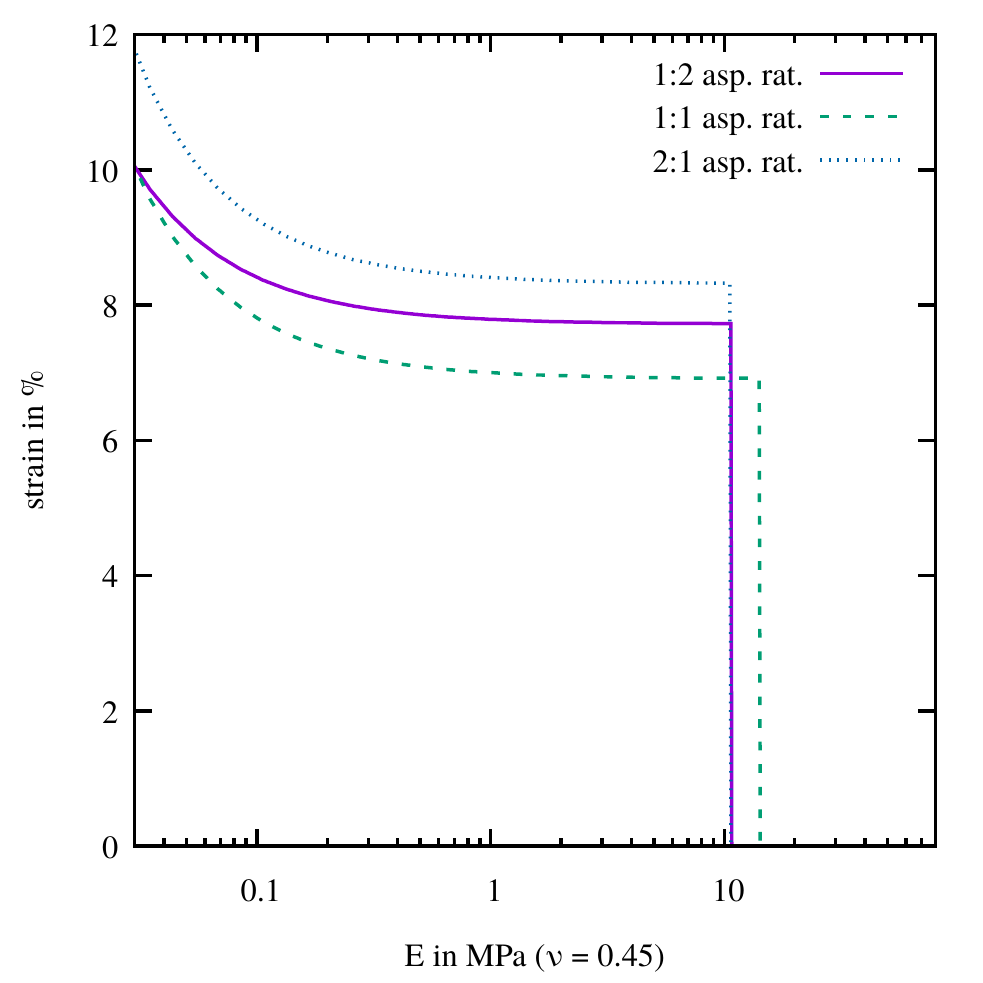}\hspace{4mm}
   \includegraphics[width=4.7cm]{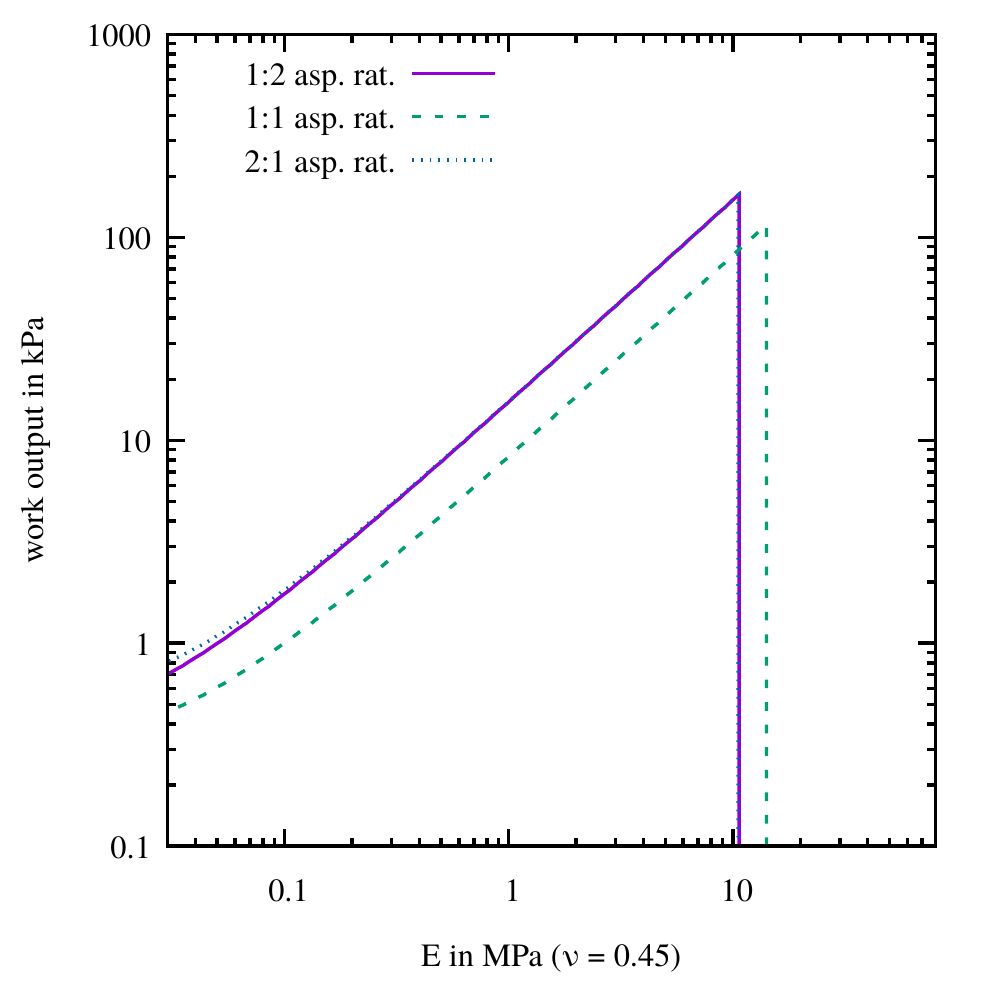}  
  \end{minipage}
 \end{center}
\caption{Macroscopic spontaneous strain and work output, as a function of the polymer's elastic modulus, for different aspect ratios of the MSM particle. The sketches on the left show the different particle shapes under consideration. The curves for 1:2 and 2:1 aspect ratios in the left panel are almost undistinguishable. The external magnetic field is horizontal.}
\label{figaspectratio}
\end{figure}

\begin{figure}
 \begin{center}
  \begin{minipage}[c]{1.3cm}
   \includegraphics[width=1.3cm]{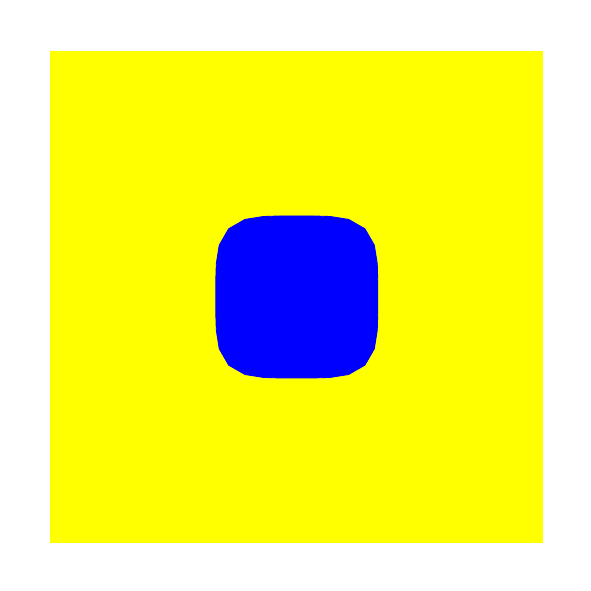}\\
   \includegraphics[width=1.3cm]{volfrac-30}\\
   \includegraphics[width=1.3cm]{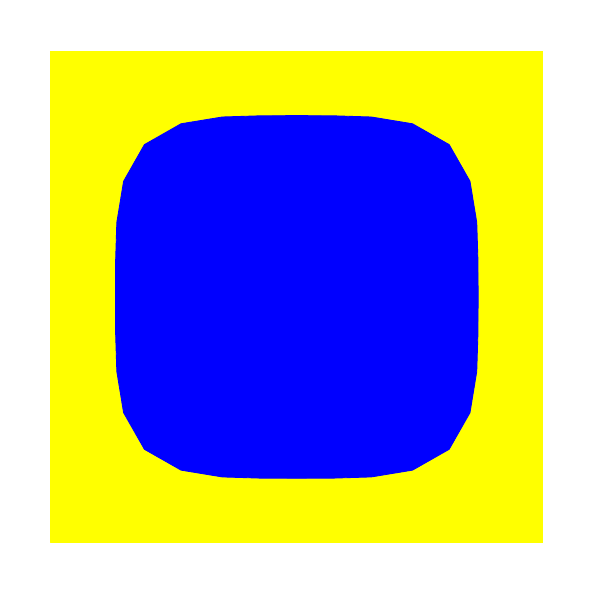}\\[0mm]
  \end{minipage}\hspace{4mm}
  \begin{minipage}[c]{9.8cm}
   \includegraphics[width=4.7cm]{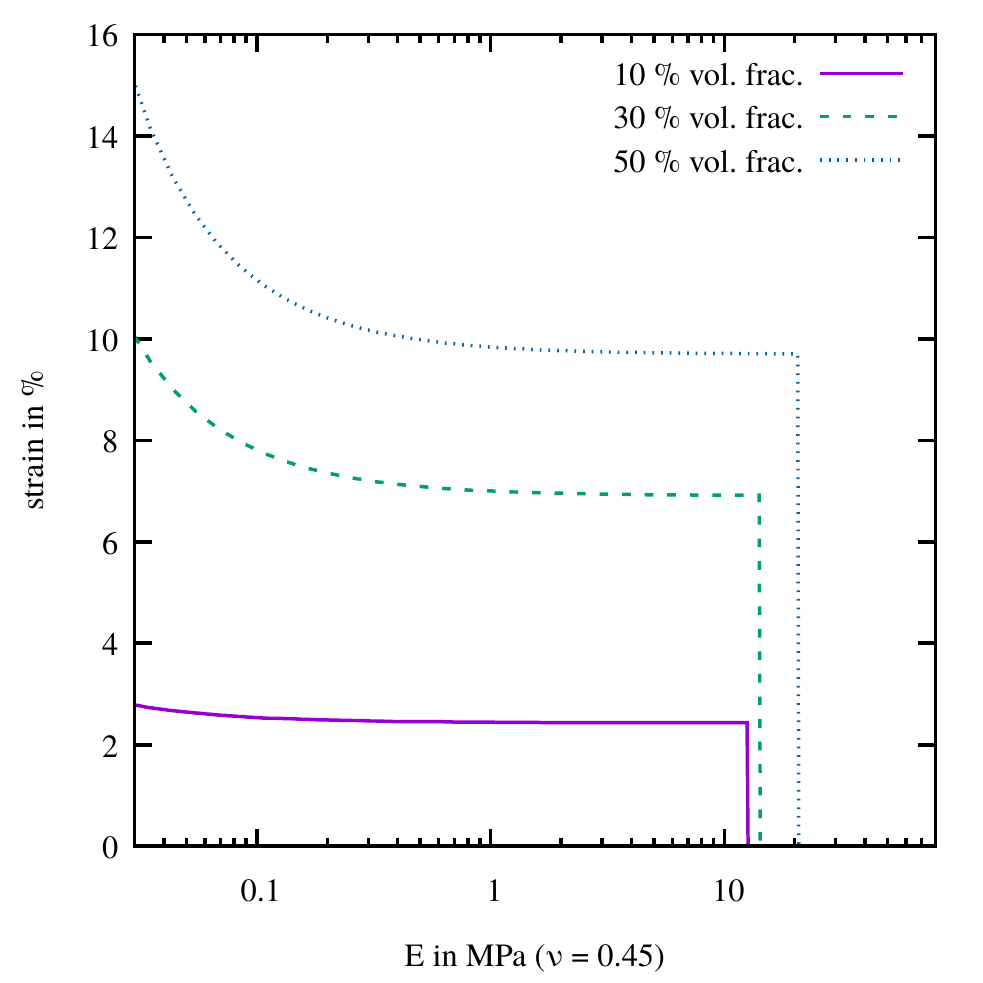}\hspace{4mm}
   \includegraphics[width=4.7cm]{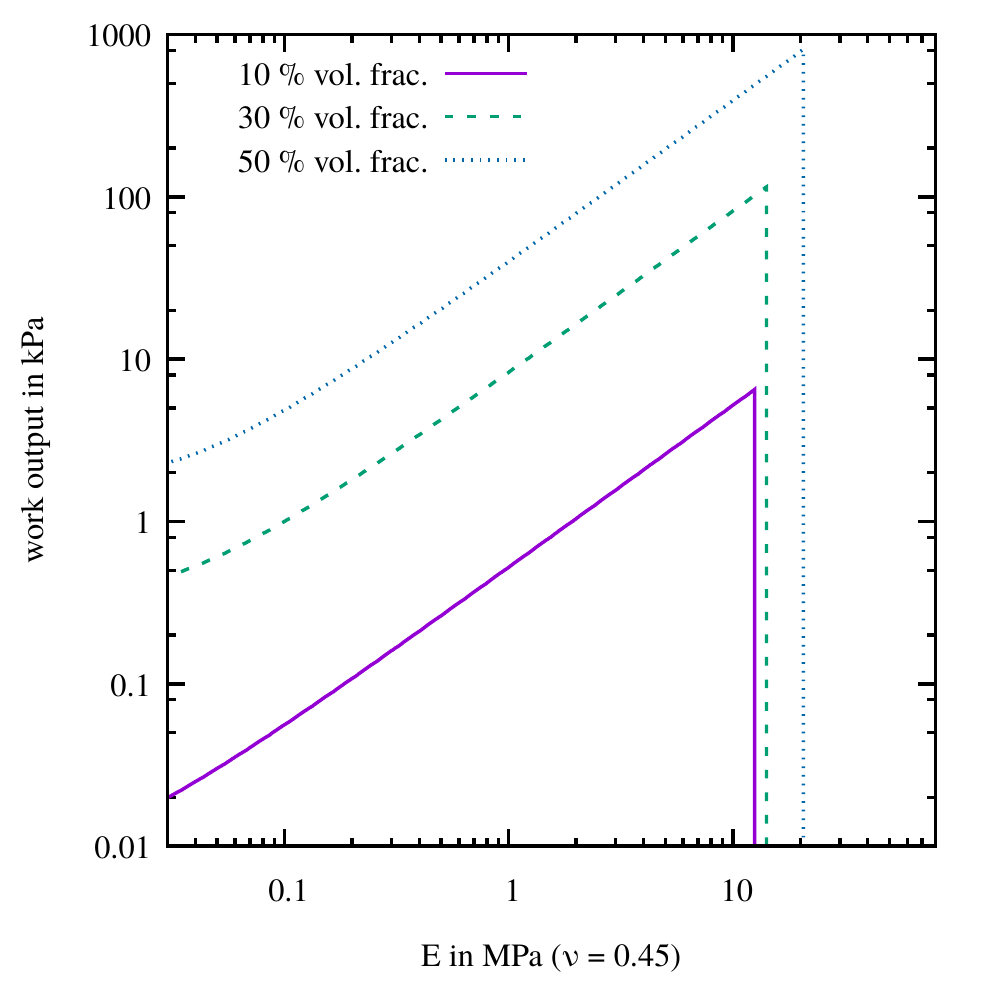}  
  \end{minipage}
 \end{center}
\caption{Macroscopic spontaneous strain and work output, as a function of the polymer's elastic modulus, for different volume fractions of MSM material. The sketches on the left show the different volume fractions.}
\label{figvolfrac}
\end{figure}

Our numerical results are illustrated in Figure \ref{figEmodul}, Figure \ref{figaspectratio} 
and Figure \ref{figvolfrac}.
We systematically display energies  and spontaneous strains, assuming affine macroscopic deformations,  in dependence on the polymer elasticity; our aim being to understand which type of polymer gives the better macroscopic material properties, for various choices of the other parameters.

The left panel of Figure  \ref{figEmodul} compares the internal energy of 
the polymer with the particles in the two different martensitic variants . The first one, called 
the ``untransformed'' phase, is the one the MSM particles had 
when the polymer was solidified, and is therefore elastically preferred. The other one is instead the one which is magnetially preferred, with the easy axis aligned with the external field.
 Comparing the two it is clear that for the parameters  considered here the phase transformation will only occur for soft polymers, with a Young modulus below around 14 MPa (black vertical line). We also plot the energy that the transformed phase would have if no macroscopic deformation would have been allowed: the difference between the two is a measure of the maximal work output that we can obtain for this material. 
The right panel illustrates the spontaneous strain corresponding to the two phases; we recall that the transformed one is only relevant for $E$ below the critical value. Both phases have a significant transformation for very soft polymers, this is a magnetostrictive effect arising from the interaction between the magnetic dipoles of the particles, and does not depend on the phase transition. It is present only for extremely soft polymer. 

Figure \ref{figaspectratio} shows the resulting spontaneous strain curve and the work output for three different shapes of the particles. Here, it is apparent that, while very soft polymers give the largest spontaneous strains, they are practically not very relevant, since the work output is minimal. The best polymers seem to be those with an intermediate Young modulus, between 1 and 10 MPa: on the one hand the transition is expected to occur reliably, since the difference between the energies of the two phases is significant, on the other hand the phase transition will lead to a macroscopic deformation, since the energy difference between the macroscopically undeformed and the macroscopically deformed state (the work output) is significant. 

Elongate particles result in a somewhat larger macroscopic strain and work output. However, the critical elastic modulus of the polymer (above which no transformation takes place) is lower for longer particles, assuming constant volume fraction. The macroscopic deformation of a composite where the particles are aligned to the magnetic field is even larger than if the particles are aligned perpendicularly, but the work output for these two configuration is nearly identical.

Finally, Figure \ref{figvolfrac} illustrates the dependence on the volume fraction. As intuitively expected, a larger amount of MSM material in the mixture leads to a larger spontaneous strain and larger work output. It is interesting to observe that the optimal choice for the polymer does not depend very strongly on the volume fraction. This is an important observation, since in material production high volume fractions are difficult to realize, and often turn out not to be uniform across the sample.


\def\no{
\section{Time-dependent extension}
\label{sectimedep}

\cite{MielkeMicroplast}

1 Bild mit 3 frames aus der neue Rechnung um das Modell zu erklären, ohne Homogenisierung

Neue Rechnungen für makro strain für unterschiedliche E-Module des Polymers
Wie Bild 11 vom 2016 Paper (Magnetisierung wechselne Richtung), aber andere vertikale axis. 
Mit Makro-Feld
}

\section{Conclusions and outlook}
In summary, we have discussed a model for composite materials, in which one component consists of magnetic-shape-memory particles, and the other of a polymer matrix. We discussed analytical homogenization, and rigorously 
derived a macroscopic effective model, whose energy density is obtained solving a suitable cell problem. We then addressed the cell problem numerically, and investigated how the microstructure can be optimized
to obtain the composite with the best material properties. 

In closing, we observe that the assumption that the particles are affinely deformed, which has a strong influence on the results presented and is only appropriate for very small particles, can be relaxed. In particular, 
in \cite{ContiLenzRumpf2012} we discussed a variant of this model in which the particles are assumed to be large with respect to the scale of the individual domains, so that three scales are present: 
the scale of a microstructure inside a single particle, the scale of the individual particles interacting with the polymer, and the scale on which macroscopic material properties are observed and measured. 
Furthermore, in \cite{ContiLenzRumpf2016} a time-dependent extension of the model was developed, 
assuming that two elastic phases are present inside each particles; the phase boundaries then move in response to the applied magnetic field. The formulation of a rate-independent model for the motion of phase boundaries permits in this case  the study of hysteresis in the composite. The much richer picture that most likely will arise in the intermediate regimes, and the extension to three spatial dimensions, still remain unexplored.

\subsection*{Acknowledgment}
This work was partially supported by the Deutsche Forschungsgemeinschaft
through Schwerpunktprogramm 1239 {\em \"Anderung von
  Mikrostruktur und Form fester Werkstoffe durch \"au{\ss}ere Magnetfelder} and through Sonderforschungsbereich 
  1060 {\em Die Mathematik der emergenten Effekte}.

\bibliographystyle{alpha-noname}
\bibliography{msm}

\end{document}